\documentclass[11pt,a4paper]{article}
\usepackage[english] {babel}
\usepackage[latin1]{inputenc}
\usepackage[T1]{fontenc}

\usepackage{amsmath, amsthm}
\usepackage{amsfonts,amssymb}
\usepackage{url}
\usepackage{graphicx}

\usepackage{latexsym}
\usepackage{bbm}
\usepackage{natbib}
\usepackage{verbatim}

\newcommand{\1}{{\mathbbm 1}}
\newcommand{\E}{{\mathbb E}}
\renewcommand{\P}{{\mathbb P}}
\newcommand{\mR}{{\mathbb R}}

\newtheorem{Rem}{Remark}

\newtheorem{Pro}{Proposition}[section]

\usepackage[margin=2.8cm]{geometry}

\begin{document}

\title{\bf Variance estimation and asymptotic confidence bands for the mean estimator of sampled functional data with high entropy unequal probability sampling designs}

\author{Herv\'e {\sc Cardot}$^{(a)}$,  Camelia {\sc Goga}$^{(a)}$ and  Pauline  {\sc Lardin}$^{(a,b)}$ \\
(a) Universit\'e de Bourgogne, Institut de Math\'ematiques de Bourgogne, \\ 9 av. Alain Savary, 21078 DIJON, FRANCE \\
(b) EDF, R{\&}D, ICAME-SOAD, 1 av. du G\'en\'eral de Gaulle, \\ 92141 CLAMART, FRANCE
}

\maketitle

\begin{abstract}
For  fixed size sampling designs with high entropy it is well known that the variance of the Horvitz-Thompson estimator can be approximated by the H\'ajek formula. The interest of this asymptotic  variance approximation is that it only involves the first order  inclusion probabilities of the statistical units. We extend this variance formula when the variable under study is functional and we prove, under general conditions on the regularity of the individual trajectories and the sampling design, that we can get a uniformly convergent estimator of the variance function of the Horvitz-Thompson estimator of the mean function. Rates of convergence to the true variance function are given for the rejective sampling.
We deduce, under conditions on the entropy of the sampling design, that it is possible to build confidence bands whose coverage is asymptotically the desired one via simulation of Gaussian processes with variance function given by the H\'ajek formula. Finally, the accuracy of the proposed variance estimator is evaluated on samples of electricity consumption data measured every half an hour over a period of one week. 
\end{abstract}

\noindent \textbf{Keywords} : covariance function, finite population, first order inclusion probabilities, H\'ajek approximation, Horvitz-Thompson estimator,   Kullback-Leibler divergence,  rejective sampling, unequal probability sampling without replacement.

\section{Introduction}

Computing the variance of the Horvitz-Thompson estimator for unequal probability sampling designs can be difficult because the variance formula involves second order probability inclusions which are not always known. 
The H\'ajek variance formula, derived in \cite{hajek_1964} for rejective sampling is an asymptotic approximation which only requires the knowledge of the first order inclusion probabilities and is easy to compute.
 It is shown in \cite{hajek_1964} and \cite{chen_al_1994} that, for given first order inclusion probabilities, the rejective sampling is the fixed size sampling design with the highest entropy. The validity of this approximation is closely related to the value of the entropy of the considered sampling design. \cite{Hajek1981} proves that this approximation is also valid for the Sampford-Durbin sampling whereas  \cite{berger_1998} gives general conditions on the relative entropy of the sampling design, also called Kullback-Leibler divergence, which justify  the use of this approximated variance formula. 
 Variants and refinements of the H\'ajek variance formula as well as variance estimators are proposed in   \cite{DevilleTille2005}. \cite{matei_tille_2005} show on simulations that these approximations to the variance of Horvitz-Thompson estimators are effective, even for moderate sample sizes, provided that the entropy of the underlying sampling design is high enough.  
Recently  \cite{DevilleTille2005} and \cite{Fuller2009} consider balanced, or approximately balanced, sampling algorithms which can be useful to build designs with fixed size and given inclusion probabilities. They relate these sampling designs to the rejective sampling, so that the H\'ajek variance approximation can be used. 
Note also that there exist other ways to get an approximation to the variance of the Horvitz-Thompson estimator which do not require the knowledge of the second order inclusion probabilities (see {\it e.g.} \cite{ShaHan2003}). These approaches do not rely on asymptotic developments  and are not considered in this work.

When the aim is to build confidence intervals, the asymptotic distribution of the Horvitz-Thompson estimator is required. 
The Central Limit Theorem has  been checked  by \cite{erdos_renyi_1959} and \cite{hajek_1960} for the  simple random sampling without replacement, by \cite{hajek_1964} for the rejective sampling and by  \cite{MR561274} for the Sampford sampling. 
\cite{berger_1998_2} states that  the Kullback Leibler divergence  of the considered sampling design, with respect to the rejective sampling, should tend to zero when the sample size gets larger for the Horvitz-Thompson estimator to be asymptotically Gaussian. 
\medskip

In recent studies in survey sampling the target was not a mean real value or a mean  vector but a mean function (see \cite{cardot_josserand_2011} and \cite{CDGJL2012} for the estimation of electricity consumption curves) and one important issue was how to build confidence bands when using $\pi$ps sampling designs. 
A rapid technique that is well adapted for large samples has been studied in \cite{Degras2011} and \cite{CardotDegrasJosserand2012}. It consists in first estimating the covariance function of the mean estimator and then simulating a Gaussian process, whose covariance function is the estimated covariance function, in order to determine the distribution of its supremum. This strategy which has been employed successfully in  \cite{CDGJL2012} to build confidence bands necessitates to have an effective estimator of the variance function of the Horvitz-Thompson estimator.
The aim of this work is to prove that under general assumptions on the sampling design and on the regularity of the trajectories, the H\'ajek formula provides a uniformly consistent estimator of the variance function.  So, it is possible to assess rigorously  confidence bands built by using the procedure described previously. 

\medskip

The paper is organized as follows. The notations and our estimators are presented in Section 2. In Section 3, we state our main result, namely the uniform convergence of the  variance function estimator obtained under broad assumptions on the regularity of the trajectories and the sampling design. We deduce that if the Horvitz-Thompson estimator of the mean curve is pointwise  asymptotically Gaussian, then it also satisfies, under the same conditions, a functional central limit theorem. The confidence bands obtained by the Gaussian process simulation techniques have  asymptotically the desired coverage. In section 4, we evaluate the performance of the covariance function estimator on samples drawn from a test population of $N=15055$ electricity consumption curves measured every half an hour over a  one-week period. 
Note there are many ways of drawing samples with high entropy sampling distribution and with given first order inclusion probabilities (see \textit{e.g.}  \cite{MR681289}, \cite{Tille2006}, \cite{MR2300911} and \cite{MR2724511}). Because of  our large population and large sample context, we use  the vast version of the Cube algorithm  (\cite{deville_tille_2004}) developed in \cite{chauvet_tille_2006} for dealing with very large populations (e.g., of millions of units). Finally, Section 6 contains some concluding remarks. 
The proofs are gathered in an Appendix.

\section{Variance estimation and the H\'ajek formula}
Let us consider a finite population $U=\{1,...,N\}$ of  known size $N$, and suppose that, for each unit $k$ of the population $U$, we can observe a deterministic curve $Y_k=(Y_k(t))_{t\in[0,T]}$.
We want  to estimate the mean trajectory $\mu_N(t)$, $t\in[0,T]$, defined as follows:
\begin{align*}
\mu_N(t) &=\frac{1}{N}\sum_{k\in U}Y_k(t).
\end{align*}

We consider a sample $s$, with fixed size $n$,  drawn from $U$ according to a sampling design $p_N(s)$, where $p_N(s)$ is the probability of drawing the sample $s$. 
The mean curve $\mu_N(t)$ is estimated by the Horvitz-Thompson estimator,
\begin{align}
\widehat{\mu}(t) & =\frac{1}{N}\sum_{k\in s}\frac{Y_k(t)}{\pi_k}=\frac{1}{N}\sum_{k\in U}\frac{Y_k(t)}{\pi_k} \1_{k}, \quad t \in [0,T],
\label{horvitz-thompson}
\end{align}
where $\1_k$ is the sample membership indicator, $\1_k=1$ if $k\in s$ and $\1_{k}=0$ otherwise. We denote by $\pi_k = \E_p(\1_k)$  the  first order inclusion probability of unit $k$ with respect to the sampling design $p_N(s)$ and we suppose that $\pi_k >0,$ for all units $k$ in $U.$ 
It is well known that, for each value of $t\in [0,T],$  $\widehat{\mu}(t)$ is a  design-unbiased estimator of $\mu_N(t)$, $i.e.$ $\E_p(\widehat{\mu}(t))=\mu_N(t).$ We denote by $\pi_{kl} = \E_p(\1_{kl})$ with $\1_{kl}=\1_k\1_l,$ the second order inclusion probabilities and we suppose that $\pi_{kl}>0$ for all $k,l\in U.$

Since the sample size is fixed, the variance $\gamma_{p}(t,t)$ for each instant $t$ of the estimator $\widehat{\mu}(t)$ is given by the Yates and Grundy formula (see \cite{YG53} and \cite{Sen53}),
 \begin{equation}
\gamma_{p}(t,t)=-\frac{1}{2}\frac{1}{N^2}\sum_{k\in U}\sum_{l\in U,l\neq k}(\pi_{kl}-\pi_k\pi_l)\left(\frac{Y_k(t)}{\pi_k}-\frac{Y_l(t)}{\pi_l}\right)^2, 
\label{variance yates-grundy0}
\end{equation}
and it is straightforward to express the covariance $\gamma_{p}(r,t)$ of  $\widehat{\mu}$ between two instants $r$ and $t,$ as follows
 \begin{equation}
\gamma_{p}(r,t)=-\frac{1}{2}\frac{1}{N^2}\sum_{k\in U}\sum_{l\in U,l\neq k}(\pi_{kl}-\pi_k\pi_l)\left(\frac{Y_k(r)}{\pi_k}-\frac{Y_l(r)}{\pi_l}\right)\left(\frac{Y_k(t)}{\pi_k}-\frac{Y_l(t)}{\pi_l}\right).
\label{variance yates-grundy}
\end{equation}

The variance formula (\ref{variance yates-grundy0}) clearly indicates that if the first order inclusion probabilities are chosen to be approximately proportional to $Y_k(t),$ the variance of the estimator $\widehat{\mu}(t)$ will be small. In practice, we can consider a non-functional auxiliary variable $X$ of values  $x_k$ supposed to be positive and  known for all the units $k\in U.$ If $X$ is nearly proportional to the  variable of interest, it can be very interesting to consider a sampling design whose  first order inclusion probabilities are given by
\begin{align*}
\pi_k &=   n\frac{x_k}{\sum_{k\in U}x_k}. 
\end{align*}
There are many ways of building sampling designs with given first order inclusion probabilities (see {\it e.g} \cite{MR681289} and \cite{Tille2006}) and we focus here on the designs with high entropy, where the entropy of a sampling design $p_N$ (a discrete probability distribution on $U$) is defined by
\begin{align*}
H(p_N) &= - \sum_{k\in s} p_N(s) \ln ( p_N(s))
\end{align*}
with the convention $0 \ln 0 = 0.$ It has been proven (see \cite{Hajek1981} and \cite{chen_al_1994}) that, for given first order inclusion probabilities, the rejective sampling, or conditional Poisson sampling, is the fixed size sampling design with the highest entropy.
Then, a key result  is the following uniform approximation to the  second order inclusion probabilities, for $k\neq l,$
\begin{align}
\pi_{kl} & =\pi_k\pi_l\left\lbrace1-\frac{(1-\pi_k)(1-\pi_l)}{d(\pi)}[1+o(1)]\right\rbrace
\label{approximation pikl}
\end{align}
where $d(\pi)=\sum_{k\in U}\pi_k(1-\pi_k)$ is supposed to tend to infinity. Note that this  implies that $n$, $N$ and $N-n$ tend to infinity. This asymptotic approximation  is satisfied  for the rejective sampling and the Sampford-Durbin sampling which is also a high entropy sampling design (see \cite{Hajek1981}). 
\begin{Rem}
Formula (\ref{approximation pikl}) can be seen rather strange and  we give an intuitive and simple interpretation  in terms of conditional covariance. Note that this is not a proof. Consider a Poisson sampling design with inclusion probabilities $p_1, \ldots, p_N$ such that $\sum_{k \in U}p_k = n$ and  $\E_p(\1_k | \#s = n) = \pi_k$ where $\# s$ denotes the (random) sample size and $\1_k$ is the indicator membership to the sample $s$ of unit $k$ (see \cite{chen_al_1994} for the existence of such sampling design).  
Considering now the  covariance given the sample size, $cov( \1_k, \1_l |  \#s = n) = \pi_{kl} - \pi_k\pi_l$, we get the following approximation, which is similar to (\ref{approximation pikl}), if we use  the formula for the conditional variance in a Gaussian framework,
\begin{align*}
cov( \1_k, \1_l |  \#s = n) &\approx cov( \1_k, \1_l ) - \frac{cov(\1_k, \# s) cov(\1_l, \# s) }{var(\# s)} \\
 &\approx 0 - \frac{\pi_k (1-\pi_k) \pi_l (1-\pi_l)}{\sum_{k\in U} \pi_k(1-\pi_k)}
\end{align*}
since $cov( \1_k, \1_l ) = 0$, $cov( \1_k, \#s) = p_k (1-p_k)$ and $var(\#s) = \sum_{k\in U} p_k (1-p_k)$ for Poisson sampling and, for each unit $k$, $p_k$ tends to $\pi_k$ as $d(\pi)$ tends to infinity (see \cite{hajek_1964}).
\end{Rem}

Then, we obtain,  for all $(r,t)\in[0,T]\times [0,T]$,  the H\'ajek approximation $\gamma_H(r,t)$ to the covariance function $cov(\hat{\mu}(t),\hat{\mu}(r))$,
by plugging in approximation (\ref{approximation pikl})  in (\ref{variance yates-grundy}), 
\begin{align}
\gamma_H(r,t) 
&=\frac{1}{N^2}\left[\sum_{k\in U}\frac{Y_k(t)Y_k(r)}{\pi_k}(1-\pi_k)-\frac{1}{d(\pi)}\left(\sum_{k\in U}(1-\pi_k)Y_k(t)\right) \left(\sum_{l\in U}(1-\pi_l)Y_l(r)\right)\right],
\end{align}
and we  consider in the following two estimators  for the covariance 
\begin{align}
\hat{\gamma}_H(r,t)&=\frac{1}{N^2}\frac{\hat{d}(\pi)}{d(\pi)}\left[\sum_{k\in s}\frac{1-\pi_k}{\pi_k^2}Y_k(t)Y_k(r)-\frac{1}{\hat{d}(\pi)}\sum_{k\in s}\left( \frac{1-\pi_k}{\pi_k}Y_k(t)\right) \sum_{l\in s}\left(\frac{1-\pi_l}{\pi_l}Y_l(r)\right)\right],
\end{align}
and $\hat{\gamma}^{*}_H(r,t) = \frac{d(\pi)}{\hat{d}(\pi)}\hat{\gamma}_H(r,t)$,
where $\hat{d}(\pi)=\sum_{k\in s}(1-\pi_k)$ is the Horvitz-Thompson estimator of $d(\pi)$. Note that $\hat{\gamma}_H(r,t)$  is  a slightly modified functional analogue of the variance estimator proposed by \cite{berger_1998}  in the real case. More exactly, the variance estimator considered by \cite{berger_1998} is $\hat{\gamma}_H(t,t)$ multiplied by  the correction factor $n/(n-1)$ so that  the  expression is exact for simple random sampling without replacement. The second estimator,  $\hat{\gamma}^{*}_H(r,t)$ is the extension to the functional case  of the \cite{DevilleTille2005}'s estimator. This latter approximation of the variance has been shown to be effective on simulation studies, even for moderate sample sizes,  by  \cite{matei_tille_2005}.

\noindent We can easily show the following property.
\begin{Pro}
If, for all $t\in[0,T]$, there is a constant $c_t$ such that $Y_k(t)=c_t\pi_k$ then $\gamma_H(r,t)=0$ and $\widehat{\gamma}_H(r,t)=\widehat{\gamma}^{*}_H(r,t)=0.$
\end{Pro}

With real data, we do not observe $Y_k(t)$ at all instants $t$ in $[0,T]$ but only for a finite set of $D$ measurement times, $0=t_1<...<t_D=T$. In functional data analysis, when the noise level is low and the grid of discretization points is fine, it is usual to perform a linear interpolation or a smoothing of the discretized trajectories in order to obtain approximations of the trajectories at every instant $t$ (see \citet{Ramsay_Silverman_Livre}). When there are nearly no measurement errors and when the trajectories are regular enough, \citet{cardot_josserand_2011} showed that  linear interpolation can provide sufficiently accurate approximations of the trajectories. Thus, for each unit $k$ in  the sample $s$, we build the interpolated trajectory 
\begin{align*}
Y_{k,d}(t) & =Y_k(t_i)+\frac{Y_k(t_{i+1})-Y_k(t_i)}{t_{i+1}-t_i}(t-t_i),\quad t\in[t_i,t_{i+1}],
\end{align*}
 and define the estimator of the mean curve $\mu_N(t)$ based on the discretized observations as follows
\begin{align}
\widehat{\mu}_{d}(t)
&= \frac{1}{N}\sum_{k\in s} \frac{Y_{k,d}(t)}{\pi_k},\quad  t\in[t_i,t_{i+1}].
\label{diffmuad}
\end{align}
Its covariance function is then estimated by
\begin{align}
\hat{\gamma}_{H,d}(r,t) &=\frac{1}{N^2}\frac{\hat{d}(\pi)}{d(\pi)}\left[\sum_{k\in s}\frac{1-\pi_k}{\pi_k^2}Y_{k,d}(t)Y_{k,d}(r)-\frac{1}{\hat{d}(\pi)}\sum_{k\in s}\left( \frac{1-\pi_k}{\pi_k} Y_{k,d}(t) \right)\sum_{l\in s}\left(\frac{1-\pi_l}{\pi_l}Y_{l,d}(r) \right)\right],
\label{gamma hat hajek d}
\end{align}
and we show in the next section that it is an uniformly consistent estimator of the variance function. Replacing $Y_k(t)$ by $Y_{k,d}(t)$ in $\hat{\gamma}^*_{H}(r,t)$, yields the variance estimator $\hat{\gamma}^*_{H,d}(r,t)$ based on the discretized values.

\section{Asymptotic properties}
All the proof are postponed in an Appendix.

\subsection{Assumptions}

 To demonstrate the asymptotic properties, we must suppose that the sample size and the population size become large. Therefore, we adopt the asymptotic approach of  \cite{hajek_1964}, assuming that $d(\pi)\rightarrow\infty$. Note that this assumption implies that $n\rightarrow\infty$ and $N-n\rightarrow\infty$. We consider a sequence of growing and nested populations $U_N$ with size $N$ tending to infinity and a sequence of samples $s_N$ of size $n_N$ drawn from $U_N$ according to the sampling design $p_N(s_N)$. The first and second order inclusion probabilities are respectively denoted by $\pi_{kN}$ and $\pi_{klN}$. For simplicity of notations and when there is no ambiguity, we drop the subscript $N$. To prove our asymptotic results we need to introduce the following assumptions.
\begin{itemize}
\item[\bf A1.] We assume that $\displaystyle\lim_{N\rightarrow \infty} \frac{n}{N}=\pi \in (0,1).$
\item[\bf A2.] We assume that $\displaystyle\min_{k \in U} \pi_k\geq\lambda>0$, $\displaystyle\min_{k \neq l\in U} \pi_{kl}\geq\lambda^*>0$ and

 \begin{align*}
\pi_{kl}& = \pi_k\pi_l\left\lbrace1-\frac{(1-\pi_k)(1-\pi_l)}{d(\pi)}[1+o(1)]\right\rbrace
\end{align*}
 uniformly in $k$ and $l.$
\item[\bf A3.] There are two positive constants  $C_2$ and $C_3$ and $\beta>1/2$ such that, for all $N$ and for all $(r,t)\in[0,T]\times [0,T]$,

$$\frac{1}{N}\sum_{k\in U}(Y_k(0))^2<C_2 \quad \mbox{and} \quad \frac{1}{N}\sum_{k\in U}(Y_k(t)-Y_k(r))^2<C_3\vert t- r\vert^{2\beta}.$$ 
\item[\bf A4.] There are two positive constants  $C_4$ and $C_5$ such that, for all $N$ and for all $(r,t)\in[0,T]\times [0,T]$,

$$\frac{1}{N}\sum_{k\in U}(Y_k(0))^4<C_4 \quad \mbox{and} \quad \frac{1}{N}\sum_{k\in U}(Y_k(t)-Y_k(r))^4<C_5 \vert t- r\vert^{4\beta}.$$ 
%
\item[{\bf A5.}] We assume that 
$$\displaystyle\lim_{N\rightarrow\infty}\max_{(k_1,l_1,k_2,l_2)\in D_{4,N}}\vert \E_p\left[(\1_{k_1l_1}-\pi_{k_1}\pi_{l_1})(\1_{k_2l_2}-\pi_{k_2}\pi_{l_2})\right]\vert = 0$$

where $\1_{kl}$ is the sample membership of the couple $(k,l)$ and $D_{4,N}$ is the set of all distinct quadruples $(i_1,...,i_4)$ from $U.$

\end{itemize}

Assumptions {\bf A1} and {\bf A2} are classical hypotheses in survey sampling and deal with the first and second order inclusion probabilities. They are satisfied for high entropy sampling designs with fixed size (see for example \cite{Hajek1981}). They directly imply that $c n \leq d(\pi) \leq n,$ for some strictly positive constant $c.$  The assumption {\bf A2} implies that  $\displaystyle\limsup_{N \rightarrow \infty}n \max_{k\neq l \in U}\vert \pi_{kl}-\pi_k\pi_l\vert<C_1<\infty.$
It also ensures that the Yates-Grundy variance estimator is always positive since $\pi_{kl} \leq \pi_k \pi_l.$

Assumption {\bf A3} and {\bf A4}  are  regularity conditions on the individual trajectories. Even if point-wise consistency, for each fixed value of $t$, can be proven without any condition on $\beta$,  these regularity conditions are required to get the uniform convergence of the mean estimator (see \cite{cardot_josserand_2011}). 
Note finally that assumption {\bf A5} is  true for SRSWOR, stratified sampling and rejective sampling (see \cite{ArratiaGoldstein2005} and \cite{BLRG2012}). More generally, it also holds  for unequal probability designs with large entropy as shown in the following proposition. Let us recall before the definition of the Kullback-Leibler divergence $K(p_N, p_{rej})$, 
\begin{align}
K(p_N, p_{rej}) &= \sum_{k\in s} p_N(s) \ln \left( \frac{p_N(s)}{p_{rej}(s)} \right), 
\end{align}
which measures how a sampling distribution $p_N(s)$  is distant from a reference sampling distribution, chosen here to be the rejective sampling $p_{rej}(s)$ since it is the design with maximum entropy for given first order inclusion probabilities. 
We can now state the following proposition which gives an upper bound of the rates of convergence to zero of the quantity in {\bf A4} in terms of Kullback-Leibler divergence with respect to the rejective sampling.  
\begin{Pro} Let $p_N$ be  a sampling design  with the same first order inclusion probabilities as $p_{rej}.$
If $d(\pi) \to \infty,$ then
\begin{align*}
\max_{(k_1,l_1,k_2,l_2)\in D_{4,N}}\left| \E_p\left[(\1_{k_1l_1}-\pi_{k_1}\pi_{l_1})(\1_{k_2l_2}-\pi_{k_2}\pi_{l_2})\right] \right| & \leq \frac{C}{d(\pi)} + \sqrt{\frac{K(p_N, p_{rej})}{2}}
\end{align*} 
for some constant $C.$
\label{preuve:A5}
\end{Pro}

A direct consequence of Proposition  \ref{preuve:A5} is that assumption  {\bf A5} is satisfied for the rejective sampling as well as for the Sampford-Durbin design, whose Kullback-Leibler divergence, with respect to the rejective sampling, tends to zero as the sample size $n$ tends to infinity (see \cite{berger_1998_2}). Note also that  the Kullback-Leibler divergence has been approximated asymptotically for other sampling designs such as the Pareto sampling in  \cite{MR2329721}.

\subsection{Convergence of the estimated variance}
Let us first recall Proposition 3.3 in \cite{cardot_josserand_2011} which states that  the estimator $\hat{\mu}_d$ is asymptotically design unbiased and uniformly convergent under mild assumptions. More precisely, if assumptions (A1)-(A3) hold and if the discretization scheme satisfies $\max_{i \in \{1,..,d_N-1\}}\vert t_{i+1}-t_i\vert^{2 \beta} =o(n^{-1})$, then for some constant $C$
\begin{align*}
\sqrt{n}\E_p\left\lbrace\sup_{t\in[0,T]}\vert\hat{\mu}_d(t)-\mu_N(t)\vert\right\rbrace &\leq C.
\end{align*}

We can now state our first result which indicates that the covariance function estimator $\hat{\gamma}_{H,d}(r,t)$ is pointwise convergent and that the variance function estimator $\hat{\gamma}_{H,d}(t,t)$ is uniformly convergent. Note that additional assumptions on the sampling design are required in order to obtain the convergence  rates.

\begin{Pro}\label{prop:gammaconv} 
\begin{enumerate}
\item Assume (A1)-(A5) hold and the sequence of discretization schemes satisfies ${\lim_{N\rightarrow\infty}\max_{i=\{1,..,d_N-1\}}\vert t_{i+1}-t_i\vert=0}$. When $N$ tends to infinity,
\begin{eqnarray}
n \E_p \left\lbrace \mid\widehat{\gamma}_{H,d}(r,t)-\gamma_{p}(r,t)\mid\right\rbrace \rightarrow 0\label{conv_gamma_1}
\end{eqnarray}
for all $(r,t) \in [0,T]\times [0,T]$ and 
\begin{eqnarray}
n \E_p \left\lbrace \sup_{t\in[0,T]}\mid\widehat{\gamma}_{H,d}(t,t)-\gamma_{p}(t,t)\mid\right\rbrace \rightarrow 0.\label{conv_gamma_2}
\end{eqnarray}
\item Under the same assumptions, the covariance function estimator $\widehat{\gamma}^{*}_{H,d}(r,t)$ satisfies (\ref{conv_gamma_1}) and the variance function estimator $\widehat{\gamma}^{*}_{H,d}(t,t)$ satisfies (\ref{conv_gamma_2}).
\end{enumerate}
\end{Pro}

A sharper result can be stated for the particular case of rejective sampling for which accurate approximations to the multiple inclusion probabilities are available (see \cite{BLRG2012}). 
\begin{Pro}
Suppose that the sample is selected with the rejective sampling design. Assume (A1)-(A4) hold and the sequence of discretization schemes satisfies \\
$\max_{i=\{1,..,d_N-1\}}\vert t_{i+1}-t_i\vert^{2\beta}=O(n^{-1}).$ Then, for all $(r,t) \in [0,T]\times [0,T]$
\begin{align*}
n^{3} \ \E_p \left[ \left(\widehat{\gamma}_{H,d}(r,t)-\gamma_{p}(r,t)\right)^2 \right] &\leq   C
\end{align*}
for some positive constant $C.$
\label{prop:speedgammaconv} 
\end{Pro}
We can note in the proof, given in the Appendix,  that the approximation error to the true variance by the H\'ajek formula is asymptotically negligible compared to the sampling error.

\subsection{Asymptotic normality and confidence bands}

Let us assume that the Horvitz-Thompson estimator of the mean curve satisfies a Central Limit Theorem for real valued quantities with new moment conditions
\begin{itemize}
\item[{\bf A6.}] There is some $\delta>0$, such that $N^{-1}\sum_{k\in U_N}\vert Y_k(t)\vert^{2+\delta}<\infty$ for all $t\in[0,T]$, and $\left\lbrace \gamma_{p}(t,t)\right\rbrace^{-1/2}\left\lbrace\hat{\mu}(t)-\mu(t)\right\rbrace\rightarrow {\cal N}(0,1)$ in distribution when $N$ tends to infinity.
\end{itemize}
This  asymptotic normality assumption is  satisfied for high entropy sampling designs (see \cite{MR561274} and \cite{berger_1998_2}).
\cite{cardot_josserand_2011} have shown that under the previous  assumptions, the central limit theorem also holds in the space of continuous functions $C[0,T]$. More precisely,  if assumptions (A1)-(A3) and (A6) hold and  the discretization points satisfy  ${\lim_{N\rightarrow\infty}\max_{i=\{1,..,d_N-1\}}\vert t_{i+1}-t_i\vert^{2 \beta} =o(n^{-1})},$ we have 
\begin{eqnarray*}
\sqrt{n}(\hat{\mu}_d-\mu)\rightarrow Z \mbox{ in distribution in }C[0,T]
\end{eqnarray*}
where $Z$ is a Gaussian random function taking values in $C[0,T]$ with mean $0$ and covariance function $\gamma_Z(r,t)=\lim_{N\rightarrow\infty}n\gamma_{p_N}(r,t).$ The reader is referred to \cite{CGL2012}  for a discussion on the reasons of using  the convergence in the space $C[0,T]. $
This important result gives a theoretical justification of the confidence bands for $\mu_N$ built as follows:
\begin{equation}\label{SCB}
\left\{ \left[ \widehat{\mu}_{d}(t) \pm c \, \frac{  \widehat{\sigma}(t)}{\sqrt{n}} \right], \: t\in [0,T] \right\},
\end{equation}
where $c$ is a suitable number and $\widehat{\sigma}(t) = \sqrt{n \widehat{\gamma}_{H,d}(t,t)}$. 

Given a confidence level $1-\alpha \in (0,1)$, one way to build such confidence bands, that is to say one way  to find an adequate value for $c_\alpha,$ is to perform simulations of centered Gaussian functions $\widehat{Z}$ defined on $[0,T]$ with mean 0 and covariance function $n \widehat{\gamma}_{H,d}(r,t)$ and then compute the quantile of order $1-\alpha$ of $\sup_{t \in [0,T]} \left| \widehat{Z}(t)/\widehat{\sigma}(t) \right|.$ In other words, we look for a cut-off point $c_\alpha$, which is random  since it depends on the estimated covariance function $\widehat{\gamma}_{H,d},$ such that
\begin{align}
\P \left(  |  \widehat{Z}(t) | \leq c_\alpha  \frac{  \widehat{\sigma}(t)}{\sqrt{n}}, \ \forall t \in [0,T] \mid \widehat{\gamma}_{H,d} \right) &= 1 - \alpha .
\end{align}
 Next proposition provides a rigorous justification for  this Monte Carlo technique which can be interpreted as parametric bootstrap:

\begin{Pro}\label{conditional weak functional convergence for Gaussian processes}
Assume (A1)-(A6) hold and the discretization scheme satisfies \\
$\max_{i=\{1,..,d_N-1\}}\vert t_{i+1}-t_i\vert^{2\beta}=o(n^{-1}).$ 

Let  $Z$ be a Gaussian process with mean zero and covariance function $\gamma_Z$. 
Let $(\widehat{Z}_N)$ be a sequence of processes such that for each $N$, conditionally on $\widehat{\gamma}_{H,d}$ defined in (\ref{gamma hat hajek d}), $\widehat{Z}_N$ is Gaussian with mean zero and covariance $n \widehat{\gamma}_{H,d}$. 
Then for all $c>0$, as $N\to\infty$, the following convergence holds in probability: 
\begin{align*}
 \mathbb{P} \left(  |\widehat{Z}_N(t)| \le c\, \widehat{\sigma} (t) ,  \: \forall t \in [0,T] \, \big| \,  \widehat{\gamma}_{H,d} \right) & \to   \mathbb{P} \left(  |Z(t)| \le c \, \sigma (t) , \: \forall t \in [0,T] \right),
 \end{align*}
 where $\widehat{\sigma}(t) = \sqrt{n \widehat{\gamma}_{H,d}(t,t)}$ and $\sigma(t) = \sqrt{\gamma_Z(t,t)}.$
\end{Pro}
The proof of Proposition \ref{conditional weak functional convergence for Gaussian processes} is very similar to the proof of Proposition 3.5  in \cite{CGL2012} and is thus omitted.
As in \cite{CardotDegrasJosserand2012}, it is possible to deduce from previous proposition that the chosen value $\widehat{c}_\alpha = c_\alpha(\widehat{\gamma}_{H,d})$ provides asymptotically the desired coverage since it satisfies 
 $$
\lim_{N\to\infty} \mathbb{P} \left( \mu(t)\in \left[ \widehat{\mu}_{d}(t) \pm \widehat{c}_\alpha \frac{\widehat{\sigma}(t)}{\sqrt{n}}  \right], \ \forall t\in [0,T] \right) = 1-\alpha.
$$


\section{Example: variance estimation for electricity consumption curves}


In this section, we evaluate the performance of the estimators $\hat{\gamma}^*_{H,d}(r,t)$ and $\hat{\gamma}_{H,d}(r,t)$ of the functional variance  $\gamma_{p}(r,t)$ of  $\hat{\mu}_d(t)$. 
Simulation studies not reported here showed that the estimators  $\hat{\gamma}^*_{H,d}(r,t)$ and $\hat{\gamma}_{H,d}(r,t)$ conduct very similarly asymptotically. This is why we only give below the simulation results for $\hat{\gamma}_{H,d}(r,t).$

We use  the same data frame as in \cite{CDGJL2012}. More exactly, we have a population $U$ of $N=15055$ electricity consumption curves measured every half an hour during one week,  so that there are $\mathcal D=336$ time points. The mean consumption during the previous week for each meter $k$, denoted $x_k$,  is used as an auxiliary variable. This variable is strongly correlated to the consumption curve $Y_k(t)$ (the pointwise correlation is always larger than 0.80) and is nearly proportional to $Y_k(t)$ at each instant $t$. It is also inexpensive to transmit.

\begin{figure}[!htbp]
\centering
\includegraphics[width=16cm]{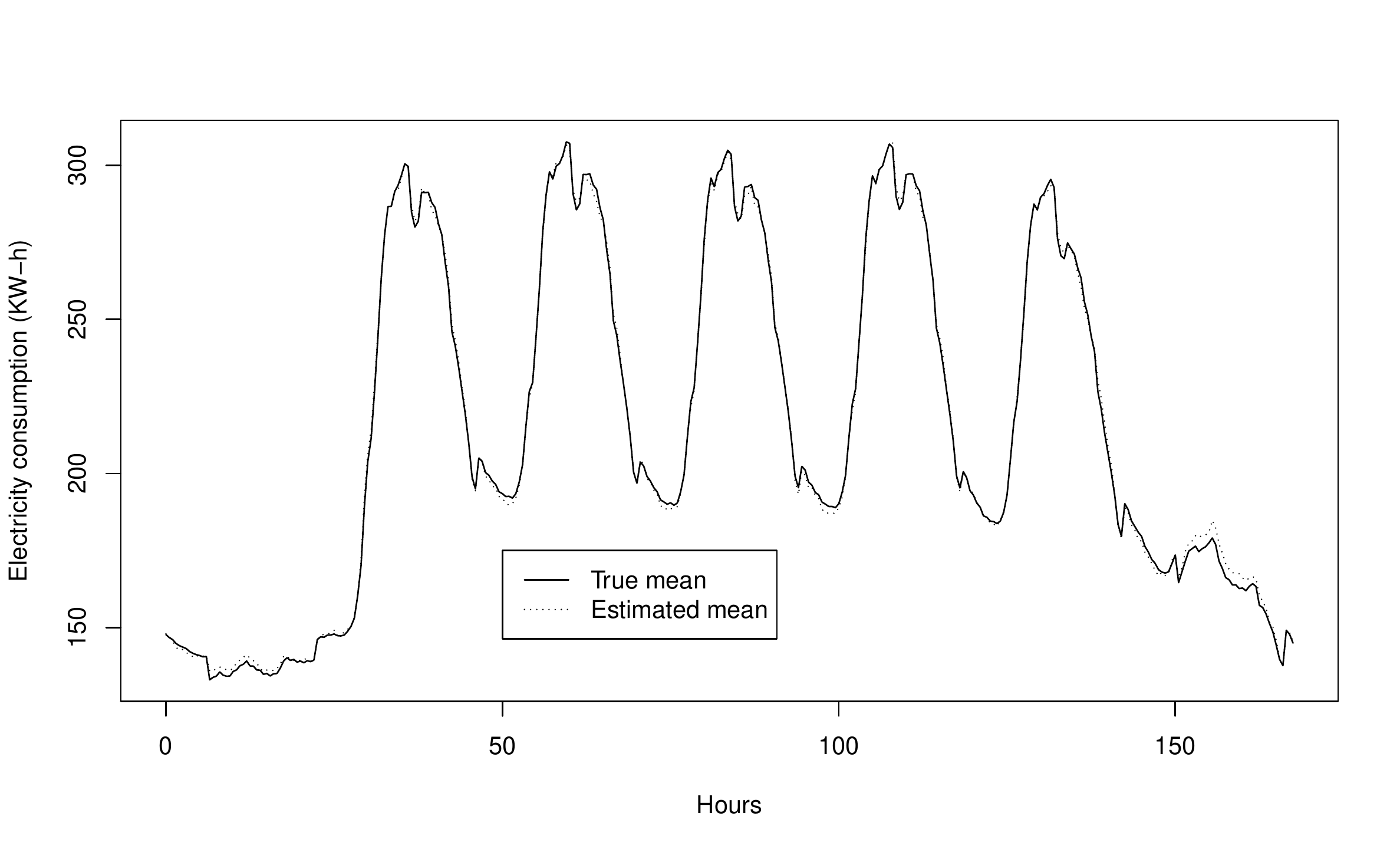}
\caption{Mean consumption curve and its Horvitz-Thompson estimation obtained from sample $s'$, with $n=1500$.}  
\label{courbe moyenne pop2}
\end{figure}

We select  samples $s$ of size $n$  drawn with inclusion probabilities $\pi_k$ proportional to the  past mean electricity consumption. This means that $\pi_k=n\frac{x_k}{\sum_{k\in U}x_k}.$ As mentioned in \cite{DevilleTille2005}, this kind of sampling may be viewed as  a balanced sampling with the balancing variable  $\boldsymbol{\pi}=(\pi_1,...,\pi_N)$. Note that by construction, the sample is also balanced on $(x_1, \ldots, x_N)$, {\em i.e} $\sum_{k\in s} x_k/\pi_k = \sum_{k\in U} x_k$.  The  sample is drawn  using the fast version (see \cite{chauvet_tille_2006}) of the cube algorithm (see \cite{deville_tille_2004}).  As suggested in \cite{chauvet_these}, a random sort of the population is made before the sample selection.  
The true mean consumption curve observed in the population $U$ and one estimation obtained from a sample $s'$ of size $n = 1500$ are drawn in Figure~\ref{courbe moyenne pop2}.



\begin{figure}[!htbp]
\centering
\includegraphics[width=16cm]{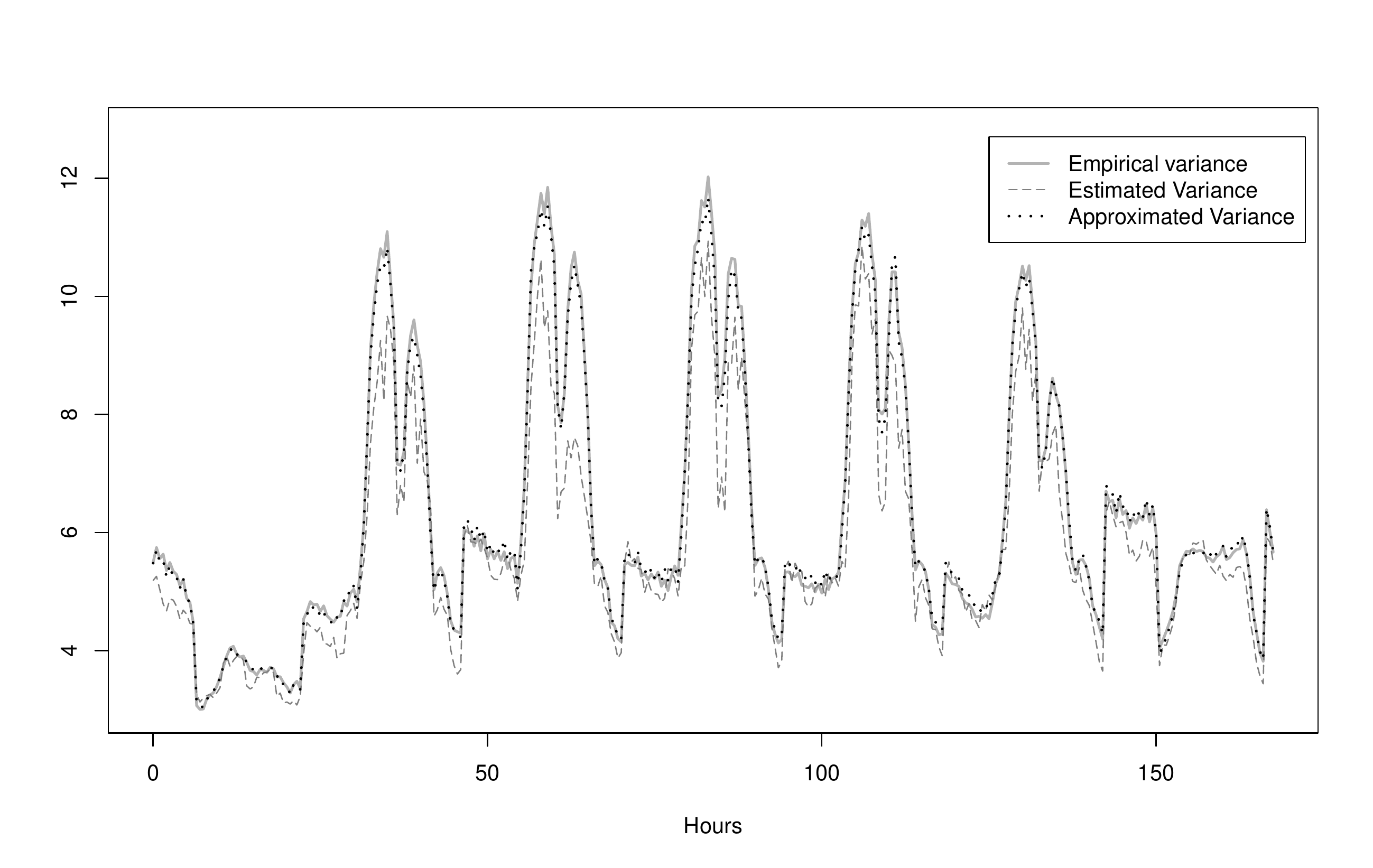}
\caption{Empirical variance $\gamma_{emp}$ (solid line), H\'ajek's approximation $\gamma_H$ (dotted line) and  variance estimation $\hat{\gamma}_{H,d}$ (dashed line) obtained from  sample $s'$, with $n=1500$. }  
\label{courbe variance hajek}
\end{figure}

The inclusion probabilities $\pi_{kl}$ being unknown, we have obtained an empirical estimation of the covariance function $\gamma_p$ via Monte Carlo. We draw $J=10000$ samples, denoted by $s_j$, for $j=1, \ldots, J$ and consider the  following Monte Carlo approximation to $\gamma_p$, 
\begin{align}
\gamma_{emp}(r,t) &=\frac{1}{J-1}\sum_{j=1}^J(\hat{\mu}_{d,j}(t)-\hat{\overline{\mu}}_d(t))(\hat{\mu}_{d,j }(r)-\hat{\overline{\mu}}_d(r)), \quad (r,t)\in[0,T]\times [0,T],
\end{align}
with $\hat{\mu}_{d,j}(t)=\frac{1}{N}\sum_{k\in s_j}\frac{Y_{k,d}(t)}{\pi_k}$, $\hat{\overline{\mu}}_d(t)=\frac{1}{J}\sum_{j=1}^J\hat{\mu}_{d,j}(t)$.
The empirical variance function $\gamma_{emp}$ (solid line) of estimator $\hat{\mu}_d$, the H\'ajek  approximation $\gamma_H$ (dotted line) and one estimation $\hat{\gamma}_{H,d}$ (dashed line) obtained from the same sample $s'$ are drawn in Figure~\ref{courbe variance hajek}.

 To evaluate the performance of  estimator $\hat{\gamma}_{H,d}$, we consider different  sample sizes, $n=250$, $n=500$ and $n=1500.$ The corresponding  values of $d(\pi)$ are  $d(\pi)=241.2,$  $d(\pi)=464.7$ and $d(\pi)=1202.3$ meaning that our asymptotic point of view is justified in this study.
 
 For each sample size, we draw $I=10000$ samples and we compute the following quadratic loss criterion
\begin{align}
R(\hat{\gamma}_{H,d})&=\frac{1}{\mathcal D}\sum_{\mathbbm{d}=1}^{\mathcal D} \frac{\vert \hat{\gamma}_H(t_\mathbbm{d},t_\mathbbm{d})-\gamma_{emp}(t_\mathbbm{d},t_\mathbbm{d})\vert^2}{\gamma_{emp}(t_\mathbbm{d},t_\mathbbm{d})^2 }\nonumber\\
&\simeq\int \frac{\vert \hat{\gamma}_H(t,t)-\gamma_{emp}(t,t)\vert^2}{\gamma_{emp}(t,t)^2 }dt.\label{erreur_R}
\end{align}

\noindent We also compute  the relative mean squared error,
\begin{align}
RMSE&=\frac{1}{I} \sum_{i=1}^IR(\hat{\gamma}^{(i)}_{H,d})\nonumber\\
&=RB^2(\hat{\gamma}_{H,d})+RV(\hat{\gamma}_{H,d}),
\end{align}
where $\hat{\gamma}^{(i)}_{H,d}$ is the value of $\hat{\gamma}_{H,d}$ computed for the $i$th simulation. It is decomposed as the sum of two terms. The term   $RB^2(\hat{\gamma}_{H,d})$ which corresponds to the square  relative bias (or approximation error) is defined by  \begin{align*}
RB(\hat{\gamma}_{H,d})^2 &=\frac{1}{\mathcal D}\sum_{\mathbbm d=1}^{\mathcal D} \left(\frac{\overline{\hat{\gamma}}_{H,d}(t_\mathbbm{d},t_\mathbbm{d})-\gamma_{emp}(t_\mathbbm{d},t_\mathbbm{d})}{\gamma_{emp}(t_\mathbbm{d},t_\mathbbm{d})}\right)^2
\end{align*}
where $\overline{\hat{\gamma}}_{H,d}(t_\mathbbm{d},t_\mathbbm{d})=\sum_{i=1}^I\hat{\gamma}^{(i)}_{H,d}(t_\mathbbm{d},t_\mathbbm{d})/I$  and $\hat{\gamma}^{(i)}_{H,d}(t_\mathbbm{d},t_\mathbbm{d})$ is the variance estimation obtained for the  $i$th simulated sample.
The second term $RV(\hat{\gamma}_{H,d}) = RMSE - RB^2(\hat{\gamma}_{H,d})$ can be interpreted as the relative variance of  estimator $\hat{\gamma}_{H,d}$.

\begin{table}[htbp]
\begin{center}
\begin{tabular}{|c|c|c|c|c|c|c|c|}
   \hline
   Sample Size & $RMSE$ & $RB^2(\hat{\gamma}_{H,d})$ & \multicolumn{5}{|c|}{$R(\hat{\gamma}_{H,d})$}\\
   \cline{4-8}
   & & &5\% &$1^{st}$ quartile& median &$3^{rd}$ quartile & 95\% \\      
     \hline\hline
    250&0.9473&0.0004&0.0188 &0.0298&0.0446&0.0748& 0.4326\\
     \hline
   500&0.3428&0.0002&0.0121&0.0191&0.0278&0.0456& 0.3510\\
 \hline
   1500&0.1406&0.0003&0.006&0.0097&0.0144&0.0272&0.0929\\
 \hline
\end{tabular}
\end{center}
\caption{$RMSE$, $RB^2(\hat{\gamma}_{H,d})$ and estimation errors according to criterion $R(\hat{\gamma}_{H,d})$ for different sample sizes, with $I=10000$ simulations. }
\label{erreur R2}
\end{table}

\begin{figure}[!htbp]
\centering
\includegraphics[width=16cm]{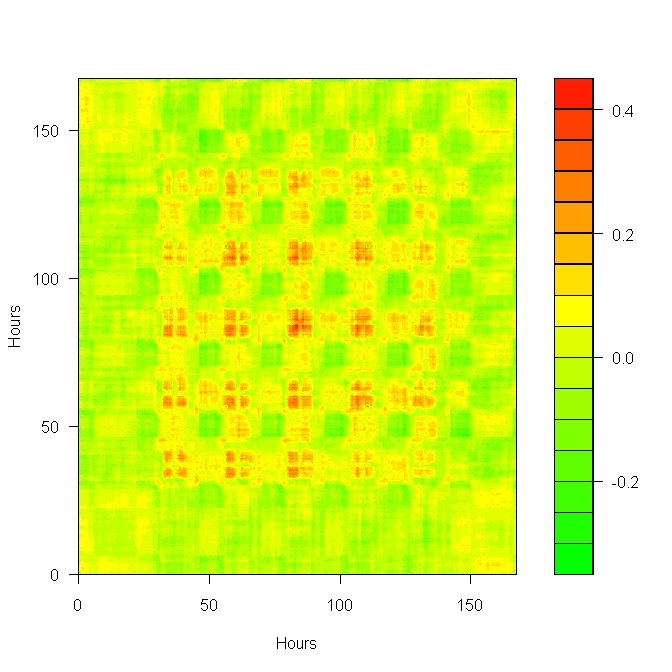}
\caption{Approximation error  $\gamma_{emp}-\gamma_{H,d}$ for a sample of size $n=1500$. }  
\label{error_emp_approx}
\end{figure}

\begin{figure}[!htbp]
\centering
\includegraphics[width=16cm]{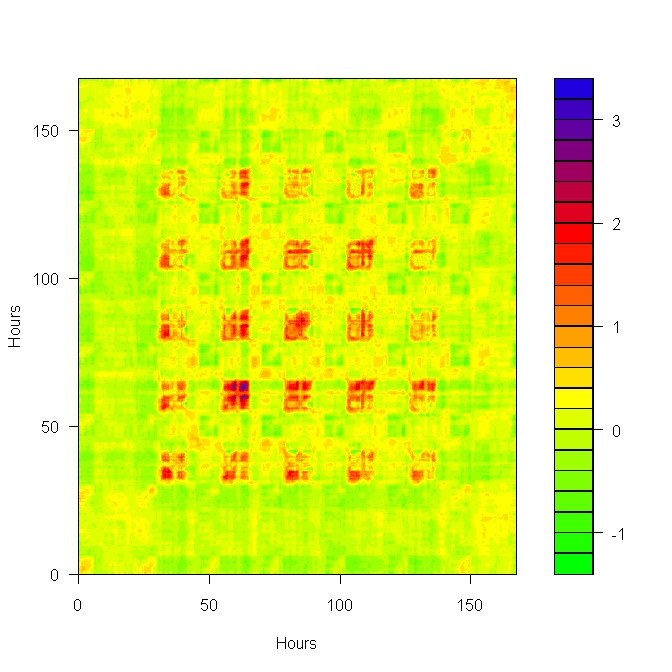}
\caption{Estimation error  $\gamma_{emp}-\hat{\gamma}_{H,d}$ for a a sample of size $n=1500$. }  
\label{error_emp_estimvarBerger}
\end{figure}

The estimation errors are presented in Table \ref{erreur R2} for the three considered sample sizes. We first note that the values of  the relative square  bias $RB(\hat{\gamma}_{H,d})$ are very low, meaning that the H\'ajek's  formula provides, in our relatively large sample context, a very good approximation to the variance. The median error for $R(\hat{\gamma}_{H,d})$ is slightly larger but remains small (always less than 5\%),  even for moderate sample sizes (n=250). This means that the most important part of the variance estimation error is due to the sampling error.   We have drawn in Figure~\ref{error_emp_approx}  the approximation error $\gamma_{emp}(t,r)-\gamma_{H,d}(t,r)$ and in Figure~\ref{error_emp_estimvarBerger} the estimation error $\gamma_{emp}(t,r)-\hat{\gamma}_{H,d}(t,r)$  for  $t,r\in \{1, \ldots, \mathcal D\},$ corresponding to a sample of size $n=1500$ with an estimation error close to the  median value of the global risk, $R(\hat{\gamma}_{H,d})= 0.0144$. It appears that the largest estimation errors for the variance occur when the level of consumption is high. We can also observe in these Figures a kind of periodic pattern which  can be related to the daily electricity consumption behavior.  

Nevertheless, we also note that the relative mean squared error $RMSE$, which is approximately equal to the relative variance of the estimator $\hat{\gamma}_{H,d},$ is rather high, especially for small sample sizes ($n=250$).  Looking at the 95 \% quantiles of $R(\hat{\gamma}_{H,d})$ in Table~\ref{erreur R2}, we can deduce that bad variance estimations only occur in rare cases but with very large errors. 
A closer look at the data shows that 
the bad performance of the variance estimator, in terms of RMSE, is  in fact due to a few individuals in the population that have both a very small  inclusion probability $\pi_k$  and a consumption level $Y_k$ that can be very important at some instants of the period. Their selection in the sample, which occurs rarely, leads to an overestimation of the mean curve and to a large error $R(\hat{\gamma}_{H,d})$ when estimating the variance at these instants.

\section{Concluding remarks}

We have studied in this work simple estimators of the covariance function of the Horvitz-Thompson estimator for curve data considering high entropy unequal probability sampling designs. Our variance estimators, which are based on  the asymptotic H\'ajek approximation to the  second order inclusion probabilities, are well suited for  large samples drawn in large populations. 
It is shown under reasonable conditions on the regularity of the curves and on the sampling design that we get  consistent estimators that  can also be used to build confidence bands for the mean, or total,  curve by employing an approach based on Gaussian process simulations. 
The illustration on the estimation of mean electricity consumption curves with $\pi$ps samples drawn with the Cube algorithm shows that, in most of  cases, the estimation error of the covariance function is small.
Nevertheless, we have in the population  a few  very influent observations (about 10 units in a population of $N=15055$) which are characterized by very small inclusion probabilities and  high values of electricity consumption at some instant of the considered period.  These influent observations, which  can be detected in the sample by considering the extreme values of the real variable $m_k = \sup_{t \in [0,T]} |Y_k(t)|/\pi_k$, completely deteriorate the quality of the variance estimator when they belong to the sample, which rarely occurs.

More robust estimators could be obtained at the sampling stage by preventing the inclusion probabilities from being too close to zero and by introducing a threshold $\delta>0$ such that 
\[
\pi_k = n \frac{\max \left(\delta,x_k \right)}{\sum_{k\in U}  \max \left(\delta,x_k \right)}.
\]
Even if the resulting Horvitz-Thompson estimator would certainly be a bit less efficient, since the proportionality would  not be respected anymore, it would permit to get a more stable estimation by attenuating the eventual effect of influent observations. 
On the other hand, another possible way to deal with this robustness issue would consist in modifying the weights of the influent observations at the estimation stage by introducing a correction such as winsorization (see \textit{e.g} \cite{MR2654641} for a review). In our variance estimation functional context, this topic is new and would certainly deserve further investigation.


\medskip

\noindent{\bf Acknowledgements}. The authors thanks the two anonymous referees as well as an associate editor for their constructive remarks.


\appendix
\section{Proofs}
Throughout the proofs we use the letter $C$ to denote a generic constant whose value may vary from place to place. Let us also define $\Delta_{kl}=\pi_{kl}-\pi_k\pi_l$ and $\Delta_{kk} = \pi_k(1-\pi_k).$ More detailed proofs can be found in \cite{Lardin2012}.

\subsection{Proof of proposition \ref{preuve:A5}}
We first consider the case of the rejective sampling $p_{rej}(s)$ and show that {\bf A5} is true if $d(\boldsymbol{\pi}_N)$ tends to infinity. 
 By Theorem 1 in in \cite{BLRG2012} and hypothesis {\bf A2}, we have 
 \begin{align*}
 \E_p \left( \1_{k_1k_2l_1l_2} \right)   -  \pi_{k_1} \pi_{k_2} \pi_{l_1} \pi_{l_2} &= O(d(\pi)^{-1})
  \end{align*}
  uniformly for $(k_1,l_1,k_2,l_2)\in D_{4,N}.$ Since $\pi_{k_1} \pi_{k_2}  - \pi_{k_1k_2} = O(d(\pi)^{-1})$ and $\pi_{l_1} \pi_{l_2} - \pi_{l_1l_2} = O(d(\pi)^{-1})$ uniformly for $(k_1,l_1,k_2,l_2)\in D_{4,N},$ we directly obtain that, for rejective sampling
  \begin{align*}
 \max_{(k_1,l_1,k_2,l_2)\in D_{4,N}}\left| \E_p\left[(\1_{k_1l_1}-\pi_{k_1}\pi_{l_1})(\1_{k_2l_2}-\pi_{k_2}\pi_{l_2})\right] \right| & \leq \frac{C}{d(\pi)}, 
  \end{align*}
  for some constant $C.$
 
 If we consider now a different sampling design $p_N(s),$ we have with
 Pinsker inequality (see  Theorem 6.1 in \cite{MR0252112}) and the property of the total variation distance,
 \begin{align*}
 \sup_{A \in {\cal A_N}} \left| p_N(A) -p_{rej}(A) \right| & \leq \sqrt{K(p_N,p_{rej})/2}
 \end{align*}
 where ${\cal A_N}$ is the set of all partitions of $U_N.$ Considering the particular cases $A = \{ (k_1,l_1,k_2,l_2)\in D_{4,N} \}$, and denoting by  $\pi_{k_1k_2l_1l_2} = p_N(A)$ and by $\pi_{k_1k_2l_1l_2}^{rej} = p_{rej}(A),$ we directly get that
 \begin{align*}
\sup_{(k_1,l_1,k_2,l_2)\in D_{4,N}}  \left| \pi_{k_1k_2l_1l_2} - \pi_{k_1k_2l_1l_2}^{rej} \right| &\leq \sqrt{K(p_N,p_{rej})/2}
 \end{align*}
and the proof is complete. 
 

\subsection{Proof of Proposition \ref{prop:gammaconv} (consistency of the covariance  and the variance functions)}

The proof follows the same steps as in \cite{CGL2012}.  We  show first that for all $t,r\in [0,T],$ the estimator of the covariance function $\widehat{\gamma}_{H,d}(r,t)$ is pointwise convergent for $\gamma_p(r,t)$  and then,  that the random variable $n(\widehat{\gamma}_{H,d}(t,t)-\gamma_{p}(t,t))$ converges  in distribution  to zero in the space $C([0,T]).$  By the definition of the convergence in distribution in $C([0,T])$ and the boundedness and continuity of the $\sup$ functional, we then directly obtain the uniformly convergence of the variance function estimator.
As in \cite{CGL2012}, in order to obtain the convergence in distribution of $n(\widehat{\gamma}_{H,d}(t,t)-\gamma_{p}(t,t)),$  we first show the convergence of all finite linear combinations which results easily from the pointwise convergence. Next, we check that the sequence $n(\widehat{\gamma}_{H,d}(t,t)-\gamma_{p}(t,t))$ is tight. 


\subsubsection*{Step 1. Pointwise convergence}

We want to show, that for each $(t,r)\in[0,T]\times [0,T]$, we have 
\begin{align*}
n \E_p\left\lbrace \mid\widehat{\gamma}_{H,d}(r,t)-\gamma_p(r,t)\mid\right\rbrace &\rightarrow  0, \quad \text{when } N\rightarrow\infty.
\end{align*}
Let us decompose
$$
n(\widehat{\gamma}_{\textrm{H},d}(r,t)- \gamma_p(r,t))= n(\widehat{\gamma}_{\textrm{H},d}(r,t)-\widehat{\gamma}_H(r,t))+n(\gamma_H(r,t)- \gamma_p(r,t))+n(\widehat{\gamma}_H(r,t)-\gamma_H(r,t))
$$
and study separately the interpolation, the approximation  and the estimation errors.

\subsubsection*{Interpolation error}

We suppose that $t\in[t_i,t_{i+1})$ and $r \in[t_{i'},t_{i'+1}).$ Using the assumptions ({\bf A1})-({\bf A3}), we can bound
\begin{align*}
n\vert\hat{\gamma}_{H,d}(r,t)-\hat{\gamma}_H(r,t)\vert& \leq C_1\vert t_{i+1}-t_i\vert^\beta+C_2\vert t_{i'+1}-t_{i'}\vert^\beta 
\end{align*}
and the assumption on the grid of discretization points leads to
\begin{align}
n\vert\hat{\gamma}_{H,d}(r,t)-\hat{\gamma}_H(r,t)\vert & =  o(1).
\label{interp:error}
\end{align}


\subsubsection*{Approximation error}
We show that, for each $(r,t) \in [0,T]\times [0,T],$
$n\mid\gamma_H(r,t)-\gamma_p(r,t)\mid  =  o(1)$. We write the approximation (\ref{approximation pikl}) as follows
\begin{align}
\pi_{kl} - \pi_k\pi_l & = - \pi_k\pi_l\frac{(1-\pi_k)(1-\pi_l)}{d(\pi)} + \frac{c_{kl}}{d(\pi)}
\label{approximationpikl2}
\end{align}
where $\max_{k\neq l\in U} | c_{kl}| \rightarrow 0$ and we use it  
in the expression of the covariance function given by (\ref{variance yates-grundy}):
\begin{align*}
\gamma_p(r,t)
&=\frac{1}{2}\frac{1}{d(\pi)N^2}\sum_{k\in U}\sum_{l\neq k\in U} \left[\pi_k\pi_l(1-\pi_k)(1-\pi_l) -c_{kl}\right]\left(\frac{Y_k(r)}{\pi_k}-\frac{Y_l(r)}{\pi_l}\right)\left(\frac{Y_k(t)}{\pi_k}-\frac{Y_l(t)}{\pi_l}\right) \\
&= \gamma_H(r,t) -\frac{1}{2}\frac{1}{N^2}\sum_{k\in U}\sum_{l\neq k\in U} \frac{c_{kl}}{d(\pi)} \left(\frac{Y_k(r)}{\pi_k}-\frac{Y_l(r)}{\pi_l}\right)\left(\frac{Y_k(t)}{\pi_k}-\frac{Y_l(t)}{\pi_l}\right). 
\end{align*}
Thus, we directly get with assumptions ({\bf A1})-({\bf A3}) that
\begin{align}
d(\pi) \  | \gamma_H(r,t)-\gamma_p(r,t) | & =    o(1).
\label{gamma h moins gamma}
\end{align}

\subsubsection*{Sampling error}
To establish the  convergence of $n(\hat{\gamma}_H(r,t)-\gamma_H(r,t))$ to zero in probability as $N\rightarrow\infty$,
it is enough to show that, for all $(r,t)\in[0,T]\times [0,T]$,
\begin{align*}
n^2\E_p\left[(\hat{\gamma}_H(r,t)-\gamma_H(r,t))^2\right] &\rightarrow 0,\quad \text{ when } N\rightarrow\infty.
\end{align*}
Noting that
\begin{align}
n\vert\hat{\gamma}_H(r,t)-\gamma_H(r,t)\vert &\leq \frac{n}{N^2}\left\vert\sum_{k\in U}\left(\frac{\hat{d}(\pi)}{d(\pi)}-1\right)\frac{\1_k}{\pi_k^2}(1-\pi_k)Y_k(t)Y_k(r)\right\vert\nonumber\\
& +\frac{n}{N^2}\left\vert\sum_{k\in U}\left(\frac{\1_k}{\pi_k}-1\right)\frac{1-\pi_k}{\pi_k}Y_k(t)Y_k(r)\right\vert\nonumber\\
& +\frac{n}{N^2}\frac{1}{d(\pi)}\left\vert\sum_{k\in U}\sum_{l\in U}\left(\frac{\1_{kl}}{\pi_k\pi_l}-1\right)(1-\pi_k)(1-\pi_l)Y_k(t)Y_l(r)\right\vert\nonumber\\
&:= |B_1(r,t)|+|B_2(r,t)|+|B_3(r,t)|,
\label{majoration gamma hat-gamma h}
\end{align}
we get
\begin{align}
n^2\E_p\left[(\hat{\gamma}_H(r,t)-\gamma_H(r,t))^2\right]&\leq 3\E_p(B_1(r,t)^2)+3\E_p(B_2(r,t)^2)+3\E_p(B_3(r,t)^2).
\label{ineq:samplinggamma}
\end{align}
Let us show now that $\E_p(B_1(r,t)^2)\rightarrow 0 $ when $N\rightarrow\infty$. Let $M=\displaystyle\max_{\pi_k\neq 1} \pi_k$.
Under the assumptions ({\bf A1})-({\bf A3}) 
and the inequality $\frac{1}{d(\pi)}\leq\frac{1}{N\lambda(1-M)}$, we have
\begin{align*}
\E_p(B_1(r,t)^2)
&\leq \frac{n^2}{\lambda^4d(\pi)^2}\E_p\Biggl[\frac{1}{N^2}(\hat{d}(\pi)-d(\pi))^2\Biggr]\left[\frac{1}{N}\sum_{k\in U} Y^2_k(t) \right] \left[\frac{1}{N}\sum_{k\in U} Y^2_k(r) \right] \leq \frac{1}{n}C
\end{align*}
since $\E_p(\frac{1}{N^2}(\hat{d}(\pi)-d(\pi))^2)=O(n^{-1}).$ Hence,  $E_p(B_1(r,t)^2)\rightarrow 0 $ when $N\rightarrow\infty$.
Now, 
\begin{align*}
\E_p(B_2(r,t)^2)&\leq \frac{n^2}{N^4}\sum_{k\in U}\sum_{l\in U}\frac{\vert\Delta_{kl}\vert}{\pi_k\pi_l}\frac{1-\pi_k}{\pi_k}\frac{1-\pi_l}{\pi_l}\vert Y_k(t)Y_k(r)Y_l(t)Y_l(r)\vert\\
&\leq \frac{1}{\lambda^3}\frac{1}{N}\left(\frac{n^2}{N^2}+\frac{n^2\max_{k\neq l\in U}\vert\Delta_{kl}\vert}{N\lambda}\right)\left(\frac{1}{N}\sum_{k\in U}\vert Y_k(t)\vert^4\right)^{1/2}\left(\frac{1}{N}\sum_{k\in U}\vert Y_k(r)\vert^4\right)^{1/2}\\
&\leq \frac{1}{N}C
\end{align*}
by assumptions (\textbf{A1})-(\textbf{A4}). Thus $\E_p(B_2(r,t)^2)\rightarrow 0 $ when $N\rightarrow\infty$.
For the third term, we have
\begin{align*}
\E_p(B_3(r,t)^2)&= n^2\E_p\Biggl[\frac{1}{N^4}\frac{1}{d(\pi)^2}\sum_{k,l\in U}\sum_{k',l'\in U}\left(\frac{\1_{kl}}{\pi_k\pi_l}-1\right)\left(\frac{\1_{k'l'}}{\pi_{k'}\pi_{l'}}-1\right)\nonumber\\
& \cdot(1-\pi_k)(1-\pi_l)(1-\pi_{k'})(1-\pi_{l'})Y_k(t)Y_l(r)Y_{k'}(t)Y_{l'}(r)\Biggr]\nonumber\\
&\leq \frac{n^2}{N^4}\frac{1}{d(\pi)^2}\sum_{k\in U}\sum_{k'\in U}\left\vert \E_p\left[\left(\frac{\1_{k}}{\pi_k^2}-1\right)\left(\frac{\1_{k'}}{\pi_{k'}^2}-1\right)\right]\right\vert\vert Y_k(t) Y_k(r)Y_{k'}(t)Y_{k'}(r)\vert\nonumber\\
& +\frac{2n^2}{N^4}\frac{1}{d(\pi)^2}\sum_{k\in U}\sum_{k'\neq l'\in U}\left\vert \E_p\left[\left(\frac{\1_{k}}{\pi_k^2}-1\right)\left(\frac{\1_{k'l'}}{\pi_{k'}\pi_{l'}}-1\right)\right]\right\vert\vert Y_k(t)Y_k(r)Y_{k'}(t)Y_{l'}(r)\vert\nonumber\\
& +\frac{n^2}{N^4}\frac{1}{d(\pi)^2}\sum_{k\neq l\in U}\sum_{k'\neq l'\in U}\left\vert \E_p\left[\left(\frac{\1_{kl}}{\pi_k\pi_l}-1\right)\left(\frac{\1_{k'l'}}{\pi_{k'}\pi_{l'}}-1\right)\right]\right\vert\vert Y_k(t)Y_l(r) Y_{k'}(t) Y_{l'}(r)\vert\nonumber\\
&:=  v_1+v_2+v_3.
\end{align*}
To bound $v_1, v_2, v_3$,  the proof follows the same lines as above. We write each double sum $\sum_{k\in U}\sum_{l\in U}$ as  the sum of two terms: the first one is $\sum_{k\in U}$ and is obtained for $k=l$ and the second one is $\sum_{k\in U}\sum_{l\neq k\in U}.$
Under assumptions ({\bf A1}), ({\bf A2}) and ({\bf A4}) and the facts that $\pi_{kl}\leq \pi_k\pi_l$ and $d(\pi)\rightarrow\infty$, we get that $v_1\rightarrow 0$. Next, we can write 
\begin{align*}
v_3
&\leq  \frac{C}{N}+\frac{n^2}{\lambda^4d^2(\pi)}\max_{(k,l,k',l')\in D_{4,N}}\left\vert \E_p\left[\left(\1_{kl}-\pi_k\pi_l\right)\left(\1_{k'l'}-\pi_{k'}\pi_{l'}\right)\right]\right\vert\left(\frac{1}{N}\sum_{k\in U}Y^2_{k}(t)\right)\left(\frac{1}{N}\sum_{l\in U}Y^2_{l}(r)\right),\nonumber 
\end{align*}
so that  $v_3\rightarrow 0$ when $N\rightarrow\infty$ and the assumptions ({\bf A1})-({\bf A5}) are fulfilled.
By the Cauchy-Schwarz inequality, we have $v_2\rightarrow 0$ when $N\rightarrow\infty .$
Finally,  we have that for all $(r,t)\in[0,T]\times [0,T]$,
$ n\vert\hat{\gamma}_H(r,t)-\gamma_H(r,t)\vert \rightarrow 0,$  when $N\rightarrow\infty$.
Finally, the proof of step 1 is complete using (\ref{interp:error}) and (\ref{gamma h moins gamma}).

\subsubsection*{Step 2. Tightness}

To check the tightness of $n(\hat {\gamma}_{H}(t,t)-\gamma_H(t,t))$  in $C[0,T]$, we use the Theorem 12.3 from  \cite{Billingsley1968} which requires that the sequence is tight for $t=0$ and that the increments of  $n(\hat {\gamma}_{H}-\gamma_H)$ between two instants $t$ and $r$ satisfy 
\begin{align*}
d^2_{\gamma}(t,r)&= n^2\E_p(\vert \widehat{\gamma}_{H}(t,t)-\gamma_H(t,t)-\widehat{\gamma}_{H}(r,r)+\gamma_H(r,r)\vert^2)  \leq C \vert t-r\vert^{2\beta}, \quad \beta>1/2
\end{align*}
for some positive constant $C$ and all $(r,t) \in [0,T]\times [0,T].$\\
The pointwise convergence of $n(\hat {\gamma}_{H}-\gamma_H)$ implies that  $n(\hat {\gamma}_{H}(0,0)-\gamma_H(0,0))$ is tight. 
Using (\ref{majoration gamma hat-gamma h}), we can decompose  $d^2_{\gamma}(t,r)$ into 3 parts,
%
%
%
\begin{align*}
d^2_\gamma(r,t)&\leq 3\Bigl(\E_p\left([B_1(t,t)-B_1(r,r)]^2\right)+\E_p\left([B_2(t,t)-B_2(r,r)]^2\right)+\E_p\left([B_3(t,t)-B_3(r,r)]^2\right)\Bigr)\\
&:= 3 \left( d^2_{B_1}+d^2_{B_2}+d^2_{B_3} \right).
\end{align*}


Denote by $\phi_{kl}(t,r) = Y_k(t)Y_l(t) -Y_k(r)Y_l(r)$ with $\phi_{k}(t,r) = Y^2_k(t) -Y^2_k(r)$ for $k=l.$  Assuming ({\bf A3}), we get that $\left(\frac{1}{N}\sum_{k\in U}|\phi_{k}(t,r)|\right)^2\leq C|t-s|^{2\beta}$ and $\left(\frac{1}{N^2}\sum_{k,l\in U}|\phi_{kl}(t,r)|\right)^2\leq C|t-s|^{2\beta}. $ Moreover, under the assumptions ({\bf A1}) and ({\bf A2}), we have 
\begin{align}
d^2_{B_1}&\leq \frac{n^2}{N^2}\left(\frac{1+\lambda}{\lambda^3}\right)^2\left(\frac{1}{N}\sum_{k\in U}|\phi_{k}(t,r)|\right)^2\leq C\vert t-r\vert^{2\beta}\label{d2 B1}
\end{align}
as well as 
\begin{align}
d^2_{B_2}&\leq \frac{n^2}{N^2}\left(\frac{1+\lambda}{\lambda^2}\right)^2\left(\frac{1}{N}\sum_{k\in U}|\phi_{k}(t,r)|\right)^2\leq C\vert t-r\vert^{2\beta}.
\label{d2 B2}
\end{align}
Finally,
\begin{align}
d^2_{B_3}
&\leq \frac{n^2}{d(\pi)^2}\left[\frac{1}{N^2}\sum_{k\in U}\sum_{l\in U}|\phi_{kl}(t,r)|\right]^2\leq C|t-r|^{2\beta}\label{d2 B3}
\end{align}
and with inequalities (\ref{d2 B1}), (\ref{d2 B2}) we deduce that
$d^2_\gamma(r,t)
\leq C\vert t-r\vert^{2\beta}.$ The proof is complete.

\vspace{0.3cm}

\noindent\textbf{Proof of Proposition \ref{prop:gammaconv}, point (2)}:
Under the assumptions ({\bf A1}) and ({\bf A2}), it is clear that $\hat d(\pi)/d(\pi)=1+o_p(1).$ 
The pointwise convergence of $n\hat{\gamma}_{H,d}^*(r,t)$ is then a direct consequence of Proposition \ref{prop:gammaconv}, point (1) and the fact that $\hat{\gamma}_{H,d}^*(r,t)=\displaystyle\frac{d(\pi)}{\hat d(\pi)}\hat{\gamma}_{H,d}(r,t)$.
Furthermore, we may write
\begin{align*}
n(\hat{\gamma}_{H,d}^*-\gamma_H)&=n\frac{d(\pi)}{\hat d(\pi)}(\hat{\gamma}_{H,d}-\gamma_H)+n\left(\frac{d(\pi)}{\hat d(\pi)}-1\right)\gamma_H.
\end{align*}
By Slutsky's theorem, the first term at  the righthand-side of previous equation converges in distribution to zero in $C([0,T])$ while the second term goes to zero in probability since $\sup_{(r,t)\in [0,T]\times [0,T]}|n\gamma_H(r,t)|<\infty$ and $\frac{d(\pi)}{\hat d(\pi)}-1=o_p(1).$  Hence, the sequence $n(\hat{\gamma}_{H,d}^*-\gamma_H)$ converges in distribution to zero in $C([0,T]). $

\subsection{Proof of Proposition \ref{prop:speedgammaconv}}

We first note that the interpolation error, bounded in (\ref{interp:error}), satisfies
\begin{align}
n^{3/2} \vert\hat{\gamma}_{H,d}(r,t)-\hat{\gamma}_H(r,t)\vert
& = O(1)
\label{vitinterp:error}
\end{align}
provided that ${\lim_{N\rightarrow\infty}\max_{i=\{1,..,d_N-1\}}\vert t_{i+1}-t_i\vert^{2\beta}=O(n^{-1})}.$
We then use the fact (see Theorem 1 in \cite{BLRG2012}) that for rejective sampling the terms $c_{kl}$ defined in (\ref{approximationpikl2}) satisfy, 
for some constant $C$,
$\max_{k,l} | c_{kl} | \leq C d(\pi)^{-1}$.
Thus, bound (\ref{gamma h moins gamma}) is now
$d(\pi)^{2} \  | \gamma_H(r,t)-\gamma_p(r,t) |  =    O(1)$.
If we examine the sampling  error, we can check that the terms $B_1$ and $B_2$ are of order $n^{-1}. $ Concerning the term $B_3,$ it is bounded by the sum $v_1+v_2+v_3$ with $v_1=O(d^{-2}(\pi))$ and $v_2\leq \sqrt{v_1v_3}.$ Thanks to Proposition \ref{preuve:A5}, we get  that the term $v_3$ satisfies $v_3 = O(d^{-1}(\pi))$ and consequently, $\E_p(B_3(r,t)^2) = O(n^{-1}).$  Thus,  
$n^2\E_p\left[(\hat{\gamma}_H(r,t)-\gamma_H(r,t))^2\right]= O(n^{-1})$
and the proof is complete.

\bibliographystyle{apalike}

\end{document}